\documentclass[12pt]{article}
\usepackage{amsfonts}
\usepackage{amssymb}
\usepackage{amsmath}

\newtheorem{theorem}{Theorem}[section]
\newtheorem{lemma}[theorem]{Lemma}

\newtheorem{rem}{Remark}

\newcommand{\qed}{\nobreak \ifvmode \relax \else
      \ifdim\lastskip<1.5em \hskip-\lastskip
      \hskip1.5em plus0em minus0.5em \fi \nobreak
      \vrule height0.75em width0.5em depth0.25em\fi}

\title{Non linear difference equations arising from a deformation of the q-Laguerre weight.  }
\author{Yang Chen\\
Department of Mathematics\\
University of Macau, (yangbrookchen@yahoo.co.uk)\\
Macau, China\\ 
James Griffin\\
Department of Mathematics\\
American University of Sharjah, (jamescgriffin@gmail.com)}

\begin{document}
\maketitle

\begin{abstract}
We study, in this paper, a one parameter deformation of the $q-$Laguerre weight function.
An investigation is made on the  polynomials orthogonal with respect to such a weight.
With the aid of the two compatibility conditions previously obtained in \cite{Chen-Ism}
 and the $q-$analog of a sum rule obtained in this paper, we derive expressions for the recurrence coefficients in terms of certain
 auxiliary quantities, and show that these quantities satisfy a pair of first order non linear difference equations.
 These difference equations are similar in form to the recognized asymmetric discrete Painleve systems such as
 $\alpha q-$P-IV and $\alpha q-$P-V.
\end{abstract}

\noindent
\section{Introduction}

In this paper we derive difference equations satisfied by the recurrence coefficients of a family of polynomials orthogonal with
 respect to the weight supported on $[0,\infty),$
\begin{equation} \label{qt-weight}
w(x,\alpha,t;q) = \frac{x^{\alpha}}{(-(1-q)x;q)_{\infty}(-(1-q)\frac{t}{x};q)_{\infty}},\quad t \geq 0, \quad \alpha > -1, \quad 0 < q < 1
\end{equation}
where
\begin{equation*}
(a;q)_{\infty} := \prod_{j=0}^{\infty}(1-a\;q^j).
\end{equation*}
In the special case of $q \to 1^-$ this reduces to
\begin{equation*}
x^{\alpha}{\rm e}^{-x}{\rm e}^{-t/x},
\end{equation*}
first considered by Chen and Its \cite{Chen-Its}, which is a singular deformation of the ordinary
 Laguerre weight. It was shown in that case that the log-derivative of the Hankel determinant (with respect to $t$) is the $\tau-$function of
 a particular P-III.
 \\
 In the limit $t \to 0^+$ our weight reduces to the $q-$Laguerre weight introduced by Moak \cite{Moak}.
 It transpires that for $t=q/(1-q)^2$ and $\alpha=0$ the corresponding orthogonal polynomials are the
 Stieltjes-Wigert polynomials. For this value of $t$ and $\alpha \neq 0$ the orthogonal polynomials were studied in \cite{Askey} by Askey.
\\
\noindent
Such a $q-$deformation causes the weight to decrease slowly near $\infty$
$$
w(x,0,t;q)={\rm O}\left({\rm e}^{-c\:(\ln\:x)^2}\right),
$$
where $c$ is positive constant,
and causes similar decrease near $0$. Weights with such slow decrease near $\infty$ were investigated in the context of unitary matrix ensembles 
which arise from electron transport in disordered systems \cite{Chen-Mutt}.
\vskip 0.2cm
\noindent
Let $\{P_n(x)\}$ be the monic polynomials orthogonal with respect to a weight $w$ on the interval $[0,\infty)$. That is
\begin{equation}
\int_0^{\infty} P_n(x)P_m(x)w(x) \; dx = h_n \; \delta_{nm},
\end{equation}
where $h_n$ is the square of the $L^2$ norm.
\vskip .15cm
\noindent
It is well known that the polynomials satisfy a three term recurrence relation,
\begin{equation}
xP_n(x) = P_{n+1}(x) + \alpha_nP_n(x) + \beta_nP_{n-1}(x).
\end{equation}
We take the initial conditions to be $\beta_0\:P_{-1}(x)=0$ and $P_0(x)=1.$
\\
Our monic polynomial has a monomial expansion
\begin{equation}
P_n(x) = x^n + p(n)x^{n-1} + ....
\end{equation}
It is clear from the recurrence relation that
\begin{equation}
\alpha_n = p(n) - p(n+1),
\end{equation}
and consequently
\begin{equation}
\sum_{j=0}^{n-1} \alpha_j = -p(n).
\end{equation}
For the weight we study, although explicit formulae are not found, we show that the recurrence
coefficients $\alpha_n$ and $\beta_n$, can be expressed in terms of the quantities $x_n$ and $y_n$ which
are solutions of a pair of coupled difference equations akin to $\alpha q-$P-IV and $\alpha q-$P-V.
In addition, we show that $p(n)$ satisfies a non linear second order difference equation in $n$.
\vskip 0.2cm
\noindent
In \cite{Chen-Ism-q}, Chen and Ismail showed that under suitable conditions on the weight $w$, the polynomials satisfied a
first order structural relation with respect to the operator $D_q,$ defined by,
\begin{equation*}
\left( D_qf\right)(x) = \frac{f(x)-f(qx)}{x(1-q)}.
\end{equation*}
Specifically, they proved the following theorem.
\begin{theorem} \label{qstruc}
Let
\begin{equation} \label{ANq}
A_n(x) = \frac{1}{h_n}\int_0^{\infty}\frac{u(qx)-u(y)}{qx-y}P_n(y)P_n(y/q)w(y) \; dy
\end{equation}
and
\begin{equation} \label{BNq}
B_n(x) = \frac{1}{h_{n-1}}\int_0^{\infty}\frac{u(qx)-u(y)}{qx-y}P_n(y)P_{n-1}(y/q)w(y) \; dy,
\end{equation}
where
\begin{equation} \label{force}
u(x) = -\frac{D_{q^{-1}}w(x)}{w(x)}.
\end{equation}
Then the orthogonal polynomials satisfy the $q-$ difference relation,
\begin{equation} \label{qstruc-eqn}
D_q P_n(x) = \beta_nA_n(x)P_{n-1}(x)-B_n(x)P_n(x).
\end{equation}
\end{theorem}
The above theorem is a $q-$analog of the structural relation appearing in \cite{Chen-Ism}.
Furthermore, it was shown in \cite{Chen-Ism-q} that the functions $A_n(x)$ and $B_n(x)$ satisfy the supplementary conditions
\begin{equation*}
B_{n+1}(x)+B_n(x) = (x-\alpha_n)A_n(x) + x(q-1)\sum_{j=0}^nA_j(x) - u(qx), \qquad \qquad  (qS_1)
\end{equation*}
and
\begin{equation*} \label{SC2}
\beta_{n+1}A_{n+1}(x) - \beta_nA_{n-1}(x) = 1+(x-\alpha_n)B_{n+1}(x) - (qx-\alpha_n)B_n(x). \qquad \qquad  (qS_2)
\end{equation*}
Equations $(qS_1)$ and $(qS_2)$ are $q-$analogs of the supplementary conditions $(S_1)$ and $(S_2)$ appearing in \cite{Chen-Ism}.
\vskip 0.2cm
\noindent
If the function $u(x)$ is rational then so are the functions $A_n(x)$ and $B_n(x)$. By comparing coefficients on both sides of
the supplementary conditions one can obtain non-linear difference equations satisfied by the recurrence coefficients.
In certain cases these equations can be solved explicitly. For example, in \cite{Chen-Ism-q} the recurrence
coefficients for the
$q$-Laguerre and Stieltjes-Wigert polynomials were extracted explicitly as a result of this
procedure. In general it is not always possible to find exact solutions, and in those cases the non linear
difference equations themselves become of interest mainly because they have been shown in many circumstances to be
related to $q$-discrete Painleve equations. For example, in \cite{Boelen-V}, Boelen and Van Assche considered the semi-classical 
$q-$Laguerre weight supported on $[0,\infty),$
\begin{equation} \label{BVA}
w(x) = \frac{x^{\alpha}(-p/x^2;q^2)_{\infty}}{(-x^2;q^2)_{\infty}(-q^2/x^2;q^2)_{\infty}}, \;\; p \in [0,q^{-\alpha}),
\; \alpha \geq 0
\end{equation}
and found that the recurrence coefficients were related to a solution of the $q-$discrete Painleve V equation. For $p=0$ the weight (\ref{BVA})
is the same as that studied by Askey in \cite{Askey} with the variable $x$ replaced by $x^2$ and the recurrence coefficients 
in this case were found to be related to a particular $q-$Painleve III. 
A non-linear difference equation
related to the following generalization of (\ref{BVA})
\begin{equation*}
w(x) = \frac{x^{\alpha}(-p_1/x^2;q^2)_{\infty}(-p_2/x^2;q^2)_{\infty}}{(-x^2;q^2)_{\infty}(-q^2/x^2;q^2)_{\infty}},
\qquad x \in [0,\infty), \;\; p_1p_2 < q^{2-\alpha}, \;\; p_1 > 0, \;\; p_2 > 0, \;\; \alpha \geq 0
\end{equation*}
was studied in \cite{Fili-Smet}.
\vskip 0.2cm
\noindent
In \cite{Ism-Mans} and \cite{Ism-John-Mans}, non-linear difference equations are derived in connection with
$q-$analogs of the Freud weights. Also in \cite{Boelen-V-Smet}, a modified $q-$Freud weight is shown to give rise to a $q-$difference
equation related to the anti-symmetric $\alpha q-$P-V equations.
\vskip 0.2cm
\noindent
In the paper of Chen and Its, \cite{Chen-Its}, they made use of the differential-difference relation,
\begin{equation}
\frac{d}{dx} P_n(x) = \beta_nA_n(x)P_{n-1}(x)-B_n(x)P_n(x).
\end{equation}
where $A_n(x)$ and $B_n(x)$ have the form
\begin{equation} \label{AN}
A_n(x) = \frac{1}{h_n}\int_0^{\infty}\frac{u(x)-u(y)}{x-y}P_n(y)P_n(y)w(y) \; dy
\end{equation}
and
\begin{equation} \label{BN}
B_n(x) = \frac{1}{h_{n-1}}\int_0^{\infty}\frac{u(x)-u(y)}{x-y}P_n(y)P_{n-1}(y)w(y) \; dy,
\end{equation}
with
\begin{equation*}
u(x) = -\frac{w^{\prime}(x)}{w(x)}.
\end{equation*}
The supplementary conditions in this case are
\begin{equation*}
B_{n+1}(x)+B_n(x) = (x-\alpha_n)A_n(x) - u(x), \qquad \qquad  \qquad \qquad \qquad (S_1)
\end{equation*}
and
\begin{equation*}
\beta_{n+1}A_{n+1}(x) - \beta_nA_{n-1}(x) = 1+(x-\alpha_n)B_{n+1}(x) - (x-\alpha_n)B_n(x). \qquad \qquad  (S_2)
\end{equation*}
Chen and Its also made use of the following equation that can be thought of as the first integral of $(S_1)$ and $(S_2)$,
\begin{equation*}
B_n^2(x) + u(x)B_n(x)+\sum_{j=0}^{n-1}A_j(x) = \beta_nA_n(x)A_{n-1}(x). \qquad \qquad \qquad (S_2^{\prime}).
\end{equation*}
A derivation of $(S_2^{\prime})$ is given in \cite{Chen-Its} and the equation first appeared in \cite{Mag}. In \cite{Chen-Its} they found a pair of coupled non-linear difference equations whose solutions were related to the recurrence coefficients. In this context they also found a particular Painleve-III differential equation in the parameter $t$. The equation $(S_2^{\prime})$ appeared in \cite{Bas-Chen}
 in connection with a Painleve V equation,
in \cite{Chen-Feig}, in connection with a Painleve IV equation and in \cite{Bas-Chen-Er},  in connection with a Painleve V equation.
\vskip 0.2cm
\noindent
For the weight appearing in (\ref{qt-weight}),  we make use of Theorem {\bf \ref{qstruc}},  as well as equations ($qS_1$) and ($qS_2$)
in an attempt to find expressions for the recurrence coefficients in terms of solutions to a pair of non-linear difference equations.
 Observe that the quantity $\sum_j A_j(x)$ appears in ($qS_1$) and not in ($qS_2$). We will find that in order to deal with
 this sum effectively we will require an additional equation involving this quantity. Therefore, instrumental in
 our approach is the derivation of a $q-$analog of the equation ($S_2^{\prime}$) which can be thought of as a first integral
 of $(qS_1)$ and $(qS_2)$. This equation appears to be new. Note that for weights such as (\ref{BVA})
 the function $u(x)$ simplifies sufficiently that the quantities involving
$\sum_j A_j(x)$ are eliminated without the use of another supplementary condition. This will not be the case for the weight
(\ref{qt-weight}).
\vskip 0.2cm
\noindent
The three main results of the paper are summarized below.
\begin{theorem} \label{LS2p}
Let $A_n(x)$ and $B_n(x)$ be given by (\ref{ANq}) and (\ref{BNq}). Then
\begin{equation*}
\beta_{n}A_{n}(x)A_{n-1}(x) = B_{n}^2(x)+u(qx)B_{n}(x) + \left(1+(1-q)xB_{n}(x)\right)\sum_{j=0}^{n-1}A_j(x). \qquad \qquad (qS_2^{\prime})
\end{equation*}
\end{theorem}
\begin{lemma} \label{rec-coeff}
Let $\{ P_n \}$ be the monic polynomials orthogonal with respect to the weight (\ref{qt-weight}) on the interval $[0, \infty)$.
Furthermore, let
\begin{equation*}
R_n = \frac{1}{h_n}\int_0^{\infty}P_n(y)P_n(y/q)\frac{w(y,\alpha,t;q)}{y} \; dy,
\end{equation*}
and
\begin{equation*}
r_n = \frac{1}{h_{n-1}}\int_0^{\infty}P_n(y)P_{n-1}(y/q)\frac{w(y,\alpha,t;q)}{y} \; dy.
\end{equation*}
Then the recurrence coefficients $\alpha_n$ and $\beta_n$ have the following form
\begin{equation*}
q^{2n+\alpha} \alpha_n =\frac{(1-q^n)}{1-q}+\frac{1-q^{n+\alpha+1}}{q(1-q)}+ q^n\left[ \frac{t}{q} \right]\left(R_n+(1-q)S_{n-1}\right),
\end{equation*}
\begin{equation*}
\beta_nq^{2n-1} = \frac{1}{q^{2\alpha}q^{2n}}\frac{1-q^n}{1-q}\frac{1-q^{n+\alpha}}{1-q} +
\frac{1-q^n}{q^{\alpha}}\left[ \frac{t}{q }\right] +
\frac{q^n}{q^{\alpha}}\left[ \frac{t}{q} \right]r_n + \frac{1}{q^{2\alpha}q^n}\left[ \frac{t}{q}\right]S_{n-1}.
\end{equation*}
where $S_{n-1}: = \sum_{j=0}^{n-1}R_j$.
\end{lemma}
\begin{rem}
The sum $S_n$ is computed in (\ref{sumR}) entirely in terms of $R_n$ and $r_n$.
Therefore the lemma above gives expressions for the recurrence coefficients in terms of $R_n$ and $r_n$ only.
\end{rem}
\begin{theorem} \label{APl}
Let
\begin{equation*}
x_n = \frac{q^{n+\alpha}(1-q)}{R_n}, \qquad \qquad y_n = q^n(1-r_n) \qquad \textrm{and} \qquad T = \frac{(1-q)^2}{q}t.
\end{equation*}
Then the $x_n$ and $y_n$ satisfy the following coupled difference equations
\begin{equation} \label{xydiff}
\begin{split}
(x_ny_n-1)(x_{n-1}y_n-1) &= q^{2n+\alpha}T\frac{(y_n-1)(y_n-1/T)}{(q^n-y_n)}, \\
(x_ny_n-1)(x_ny_{n+1}-1) &= -q^{2n+\alpha+1}\frac{(x_n-1)(x_n-T)}{x_n}.
\end{split}
\end{equation}
\end{theorem}
Note that these are similar in form to  $\alpha q-$P-IV, \cite{Van-A}, and $\alpha q-$P-V, \cite{Boelen-V-Smet}, which we state below:
\begin{equation*}
\begin{split}
\alpha q-PIV \qquad (x_ny_n-1)(x_{n-1}y_n-1) &= \frac{(y_n-a)(y_n-b)(y_n-c)(y_n-d)}{(y_n-\kappa \rho_n)(y_n-\rho_n/\kappa)} \\
(x_ny_n-1)(x_ny_{n+1}-1) &= \frac{(x_n-1/a)(x_n-1/b)(x_n-1/c)(x_n-1/d)}{(x_n-\mu w_n)(x_n-w_n/\mu)}
\end{split}
\end{equation*}
\vskip 0.2cm
\noindent
\begin{equation*}
\begin{split}
\alpha q-PV \qquad (x_ny_n-1)(x_{n-1}y_n-1) &= q^{2n}\frac{(y_n-a)(y_n-b)(y_n-c)(y_n-d)}{(q^n-\kappa y_n)(q^n-y_n/\kappa)} \\
(x_ny_n-1)(x_ny_{n+1}-1) &= q^{2n+1}\frac{(x_n-1/a)(x_n-1/b)(x_n-1/c)(x_n-1/d)}{(q^{n+1/2}-\mu y_n)(q^{n+1/2}-y_n/\mu)}.
\end{split}
\end{equation*}
Here $a,b,c,d,\kappa$ and $\mu$ are parameters, and $\rho_n$ and $w_n$ are in \cite{greeks}. Furthermore, $abcd=1.$ 
\vskip 0.2cm
\noindent
This paper is organized as follows. In Section $2$ we evaluate the rational functions $A_n$ and $B_n$ in terms of certain auxiliary quantities.
In section $3$ we give a proof of lemma {\bf \ref{LS2p}} and derive expressions for the recurrence coefficients in terms these auxiliary quantities. In section $4$ we give a proof of theorem {\bf \ref{APl}}.
Finally in section $5$ we derive a second order, non-linear difference equation, for the quantity $p(n)$.

\section{The Structural Relation}
\setcounter{equation}{0}

In this section we compute the functions $A_n(x)$ and $B_n(x)$ appearing in the relation (\ref{qstruc-eqn}) in terms of
certain auxiliary quantities. We will make use of the $q-$product rule and $q-$integration by parts, namely,
\begin{equation*}
D_q(f(x)g(x)) = f(qx)D_qg(x) + g(x)D_qf(x),
\end{equation*}
\vskip 0.4cm
\noindent
\begin{equation*}
\int_0^{\infty}f(x)D_qg(x) \; dx = -\frac{1}{q}\int_0^{\infty}g(x)D_{q^{-1}}f(x) \; dx.
\end{equation*}
The second formula is valid whenever the integrals
\begin{equation*}
\int_0^{\infty}f(x)g(x)\frac{dx}{x} \qquad \textrm{and} \qquad \int_0^{\infty}f(x)g(qx)\frac{dx}{x}.
\end{equation*}
exist.
Our first task is to compute the function $u$ for the weight function (\ref{qt-weight}). We have,
\begin{equation*}
u(x,\alpha,t;q) = \frac{x^2 + \left[\frac{1-q^{-\alpha}}{1-q}\right]qx-q^{1-\alpha}t}{x^2(1+(q^{-1}-1)x)}.
\end{equation*}
From this, it follows that,
\begin{equation}
\begin{split} \label{final}
\frac{u(qx,\alpha,t;q)-u(y,\alpha,t;q)}{qx-y} &= \frac{1}{q^{\alpha}}\left[\frac{t}{qy}\right]\frac{1}{x^2} +
\frac{1}{q^{\alpha}}\left[ \frac{q-t(1-q)^2}{q^2(1+(q^{-1}-1)y)}\right]
\frac{1}{x} \\
& \qquad -\frac{1}{q^{\alpha}}\left[ \frac{q-t(1-q)^2}{q^2(1+(q^{-1}-1)y)}\right]\frac{1-q}{1+(1-q)x}-
\left[ \frac{u(y,\alpha,t;q)}{q}\right]\frac{1}{x}.
\end{split}
\end{equation}
Note that (\ref{final}) is a rational function of $x$ and $y$.
Now we look at the effect of the $u(y,\alpha,t;q)$ term in (\ref{final}) on the quantities $A_n(x)$ and $B_n(x)$.
First we start with the function $A_n(x)$. For the $u(y,\alpha,t;q),$ term we have
\begin{equation*}
\begin{split}
\frac{1}{h_n} \int_0^{\infty}u(y,\alpha,t;q)P_n(y)P_n(y/q)w(y,\alpha,t;q) &= -\frac{1}{h_n}\int_0^{\infty} P_n(y)P_n(y/q)
D_{q^{-1}}w(y,\alpha,t;q) \; dy \\
&= \frac{q}{h_n}\int_0^{\infty}D_q \left[ P_n(y)P_n(y/q) \right] w(y,\alpha,t;q) \; dy \\
&= \frac{q}{h_n}\int_0^{\infty}\left[ P_n(y)D_q P_n(y)  + P_n(y)D_q P_n(y/q) \right] w(y,\alpha,t;q) \\
& = 0
\end{split}
\end{equation*}
where the last line follows from orthogonality. We repeat the calculation for the $u(y,\alpha,t;q)$ term appearing in  $B_n(x)$, and obtain,
\begin{equation*}
\begin{split}
& \frac{1}{h_{n-1}}\int_0^{\infty}u(y,\alpha,t;q)P_n(y)P_{n-1}(y/q)w(y,\alpha,t;q)  \\
& \qquad \qquad \qquad = -\int_0^{\infty} P_n(y)P_{n-1}(y/q)D_{q^{-1}}w(y,\alpha,t;q) \; dy \\
& \qquad \qquad \qquad = \frac{q}{h_{n-1}}\int_0^{\infty}D_q \left[ P_n(y)P_{n-1}(y/q) \right] w(y,\alpha,t;q) \; dy \\
& \qquad \qquad \qquad = \frac{q}{h_{n-1}}\int_0^{\infty}\left[ P_{n-1}(y)D_q P_n(y)  + P_n(y)D_q P_{n-1}(y/q) \right] w(y,\alpha,t;q) \\
& \qquad \qquad \qquad = q\frac{1-q^n}{1-q}.
\end{split}
\end{equation*}
With the following definitions of the auxiliary quantities
\begin{equation*}
R^{(1)}_n = \frac{1}{h_n}\int_0^{\infty}P_n(y)P_n(y/q)\frac{w(y,\alpha,t;q)}{y} \; dy
\end{equation*}
\vskip 0.2cm
\noindent
\begin{equation*}
R^{(2)}_n = \frac{1}{h_n}\int_0^{\infty}P_n(y)P_n(y/q)\frac{w(y,\alpha,t;q)}{1+y(q^{-1}-1)} \; dy
\end{equation*}
\vskip 0.2cm
\noindent
\begin{equation*}
r^{(1)}_n = \frac{1}{h_{n-1}}\int_0^{\infty}P_n(y)P_{n-1}(y/q)\frac{w(y,\alpha,t;q)}{y} \; dy
\end{equation*}
\vskip 0.2cm
\noindent
\begin{equation*}
r^{(2)}_n = \frac{1}{h_{n-1}}\int_0^{\infty}P_n(y)P_{n-1}(y/q)\frac{w(y,\alpha,t;q)}{1+y(q^{-1}-1)} \; dy
\end{equation*}
we find that $A_n(x)$ and $B_n(x)$ appearing in (\ref{qstruc-eqn}) are rational functions of $x$, and read,
\begin{equation} \label{A}
A_n(x) = \frac{R^{(1)}_n}{x^2}\left[ \frac{t}{q} \right]\frac{1}{q^{\alpha}} +\frac{R^{(2)}_n}{x}
\left[ \frac{q-t(1-q)^2}{q^2}\right]\frac{1}{q^{\alpha}} - (1-q)\frac{R^{(2)}_n}{1+x(1-q)}\left[ \frac{q-t(1-q)^2}{q^2}\right]\frac{1}{q^{\alpha}}
\end{equation}
\vskip 0.2cm
\noindent
\begin{equation} \label{B}
B_n(x) = \frac{r^{(1)}_n}{x^2}\left[ \frac{t}{q} \right]\frac{1}{q^{\alpha}} +\frac{r^{(2)}_n}{x}\left[ \frac{q-t(1-q)^2}{q^2}\right]\frac{1}{q^{\alpha}}
- (1-q)\frac{r^{(2)}_n}{1+x(1-q)}\left[ \frac{q-t(1-q)^2}{q^2}\right]\frac{1}{q^{\alpha}} - \frac{1}{x}\left[ \frac{1-q^n}{1-q} \right].
\end{equation}
We now take note of the fact that the $R_n^{(2)}$ can be expressed in terms of $R_n^{(1)}$.
Likewise for $r_n^{(2)}$ in terms of $r_n^{(1)}$.
To see this, observe that,
\begin{equation}
q^{-\alpha}\left( (1-q)\left[ \frac{t}{q}\right]\frac{1}{y}+ \left[ \frac{q-t(1-q)^2}{q^2}\right]\frac{1}{1+(q^{-1}-1)y}
   \right) w(y,\alpha,t;q) = \frac{1}{q}w(y/q,\alpha,t;q).
\end{equation}
Therefore we have
\begin{equation} \label{switchR}
q^{-\alpha}\left( (1-q)\left[ \frac{t}{q}\right]R_n^{(1)}+ \left[ \frac{q-t(1-q)^2}{q^2}\right]R_n^{(2)}  \right) = q^n.
\end{equation}
\begin{equation} \label{switchr}
q^{-\alpha}\left( (1-q)\left[ \frac{t}{q}\right]r_n^{(1)}+ \left[ \frac{q-t(1-q)^2}{q^2}\right]r_n^{(2)}  \right) = -(1-q)q^{n-1}\sum_{j=0}^{n-1}\alpha_j.
\end{equation}
\begin{rem} \label{Askey-coef}
If we set $t=q/(1-q^2)$ as was the case in \cite{Askey} then the coefficients of $R_n^{(2)}$ and $r_n^{(2)}$ in
(\ref{switchR}) and (\ref{switchr}) vanish, leading immediately to an explicit expression for $R_n^{(1)}$.
Substituting this value into the first equation of lemma {\bf \ref{rec-coeff}} gives an explicit expression for $\alpha_n$.
Subsequently, substituting this expression into (\ref{switchr}) gives an explicit expression for $r_n^{(1)}$ which can then be
substituted into the second equation in lemma {\bf \ref{rec-coeff}} to give an explicit expression for $\beta_n$.
\end{rem}
\vskip 0.2cm
\noindent
Using (\ref{switchR}) and (\ref{switchr}) we now eliminate $R^{(2)}_n$ and $r^{(2)}_n$ from (\ref{A}) and (\ref{B}). Consequently,
\begin{equation} \label{A2}
\begin{split}
A_n(x) &= \frac{R^{(1)}_n}{x^2}\left[ \frac{t}{q} \right]\frac{1}{q^{\alpha}}
+\frac{1}{x}\left(q^n - \frac{1}{q^{\alpha}}(1-q)\left[ \frac{t}{q}\right]R_n^{(1)} \right) \\
& \qquad \qquad \qquad -\frac{1-q}{1+x(1-q)}\left(q^n - \frac{1}{q^{\alpha}}(1-q)\left[ \frac{t}{q}\right]R_n^{(1)} \right),
\end{split}
\end{equation}
and
\begin{equation} \label{B2}
\begin{split}
B_n(x) &= \frac{r^{(1)}_n}{x^2}\left[ \frac{t}{q} \right]\frac{1}{q^{\alpha}}
+\frac{1}{x}\left(-(1-q)q^{n-1}\sum_{j=0}^{n-1}\alpha_j - \frac{1}{q^{\alpha}}(1-q)\left[ \frac{t}{q}\right]r_n^{(1)} \right) \\
& \qquad \qquad -\frac{1-q}{1+x(1-q)}\left(-(1-q)q^{n-1}\sum_{j=0}^{n-1}\alpha_j
- \frac{1}{q^{\alpha}}(1-q)\left[ \frac{t}{q}\right]r_n^{(1)} \right) \\
& \qquad \qquad \qquad - \frac{1}{x}\left[ \frac{1-q^n}{1-q} \right].
\end{split}
\end{equation}

\section{The Recurrence Coefficients}
\setcounter{equation}{0}
In this section we derive expressions for the recurrence coefficients in terms of the quantities
$R^{(1)}_n$ and $r^{(1)}_n$. Because we no longer require 
$R^{(2)}_n$ and $r^{(2)}_n$ we will drop the superscript and use the notation
\begin{equation*}
R_n = R^{(1)}_n, \qquad r_n = r^{(1)}_n \qquad \textrm{and} \qquad S_n = \sum_{j=0}^n R_n.
\end{equation*}
We begin with the derivation of ($qS_2^{\prime}$).
\vskip 0.2cm
\noindent
{\bf Proof of Theorem \ref{LS2p}}
\vskip 0.2cm
\noindent
First we write ($qS_2$) in the form
\begin{equation*}
\beta_{n+1}A_{n+1}(x)-\beta_nA_{n-1}(x) = 1+(x-\alpha_n)(B_{n+1}(x)-B_n(x))+(1-q)xB_n(x).
\end{equation*}
If we multiply the above equations by $A_n(x)$ and use ($qS_1$) to substitute for $(x-\alpha_n)A_n(x),$ we obtain
\begin{equation*}
\begin{split}
& \beta_{n+1}A_{n+1}(x)A_n(x) - \beta_nA_n(x)A_{n-1}(x)  = \\
& \qquad \qquad A_n(x) + (B_{n+1}^2(x)+u(qx)B_{n+1}(x)) - (B_{n}^2(x)+u(qx)B_{n}(x))\\
& \qquad \qquad \qquad + x(1-q)\left(B_{n+1}(x)\sum_{j=0}^nA_j(x) - B_n(x) \sum_{j=0}^{n-1}A_j(x) \right).
\end{split}
\end{equation*}
Observe that, up to $A_n(x)$ on the right side, the above is a first order difference equation in $n$, hence, summing over $n$,
we obtain the $q$-analog of $(S_2^{\prime})$
\begin{equation*}
\beta_{n}A_{n}(x)A_{n-1}(x) = B_{n}^2(x)+u(qx)B_{n}(x) + \left(1+(1-q)xB_{n}(x)\right)\sum_{j=0}^{n-1}A_j(x).
\end{equation*}
\hfill $\square$
\vskip 0.2cm
\noindent
To proceed further, we obtain, equating the coefficients of $x^{-2}$ in $(qS_1)$, :
\begin{equation} \label{one}
r_{n+1} + r_n = -\alpha_nR_n +1.
\end{equation}
Equating the coefficients of $x^{-1}$ in $(qS_1),$ we obtain the equation :
\begin{equation*}
\begin{split}
&-(1-q)\left(q^n \sum_{j=0}^n \alpha_j + q^{n-1}\sum_{j=0}^{n-1}\alpha_j  \right)-\frac{1}{q^{\alpha}}(1-q)\left[ \frac{t}{q} \right]\left( r_{n+1} + r_n \right) - \frac{1-q^{n+1}}{1-q} - \frac{1-q^n}{1-q} \\
& \qquad = -\alpha_nq^n - \frac{1-q^{-\alpha}}{1-q} + \frac{1}{q^{\alpha}}\left[ \frac{t}{q} \right]\left(R_n + \alpha_n(1-q)R_n-(1-q)S_n - (1-q) \right).
\end{split}
\end{equation*}
The above can be simplified making use of (\ref{one}), and we arrive at
\begin{equation*}
\begin{split}
&-(1-q)\left(q^n \sum_{j=0}^n \alpha_j + q^{n-1}\sum_{j=0}^{n-1}\alpha_j  \right) - \frac{1-q^{n+1}}{1-q} - \frac{1-q^n}{1-q} \\
& \qquad \qquad \qquad = -\alpha_nq^n - \frac{1-q^{-\alpha}}{1-q} + \frac{1}{q^{\alpha}}\left[ \frac{t}{q} \right]\left(R_n -(1-q)S_n \right).
\end{split}
\end{equation*}
The above equation simplifies to
\begin{equation*}
q^{n+1}\sum_{j=0}^n \alpha_j - q^{n-1}\sum_{j=0}^{n-1}\alpha_j = \frac{1-q^{n+1}}{1-q} + \frac{1-q^n}{1-q} - \frac{1-q^{-\alpha}}{1-q}+ \frac{1}{q^{\alpha}}\left[ \frac{t}{q} \right]\left(qS_n - S_{n-1} \right)
\end{equation*}
We multiply both sides of this equation by the integrating factor $q^{n-1}$ to obtain
\begin{equation*}
q^{2n}\sum_{j=0}^n \alpha_j - q^{2n-2}\sum_{j=0}^{n-1}\alpha_j = q^{n-1}\frac{1-q^{n+1}}{1-q} + q^{n-1}\frac{1-q^n}{1-q} - q^{n-1}\frac{1-q^{-\alpha}}{1-q}+ \frac{1}{q^{\alpha}}\left[ \frac{t}{q} \right]\left(q^nS_n - q^{n-1}S_{n-1} \right).
\end{equation*}
Summing this equation we obtain
\begin{equation} \label{pcoef}
q^{2n}\sum_{j=0}^n \alpha_j = \frac{1}{q}\left( \frac{1-q^{n+1}}{1-q}\right)^2 - \frac{1}{q}\left( \frac{1-q^{n+1}}{1-q}\right)\left( \frac{1-q^{-\alpha}}{1-q} \right)+ \frac{1}{q^{\alpha}}\left[ \frac{t}{q} \right]q^nS_n.
\end{equation}
\begin{rem}
Because $\sum_{j=0}^n \alpha_j = -p(n+1)$, (\ref{pcoef}) gives the summation $\sum_{j=0}^n R_j$ in a closed form.
\end{rem}
From (\ref{pcoef}), we see that,
\begin{equation*}
q^{2n} \alpha_n = 2\left(\frac{1-q^n}{1-q}\right) + \frac{1}{q} - \left(\frac{1}{q}+1-q^n\right)\left( \frac{1-q^{-\alpha}}{1-q} \right)+ \frac{1}{q^{\alpha}}\left[ \frac{t}{q} \right]q^n\left(S_n-qS_{n-1}\right)
\end{equation*}
which we write as
\begin{equation} \label{two}
q^{2n} \alpha_n = 2\left(\frac{1-q^n}{1-q}\right) + \frac{1}{q} - \left(\frac{1}{q}+1-q^n\right)\left( \frac{1-q^{-\alpha}}{1-q} \right)+ \frac{1}{q^{\alpha}}\left[ \frac{t}{q} \right]q^n\left(R_n+(1-q)S_{n-1}\right).
\end{equation}
Equating the coefficients of of $(1+x(1-q))^{-1}$ in $(qS_1)$ we obtain (\ref{two}) again. We go through the same process for $(qS_2^{\prime})$ and expect to find three more equations. Equating the coefficients of $x^{-4}$ in $(qS_2^{\prime})$ gives
\begin{equation} \label{three}
\beta_nR_nR_{n-1} = r_n^2 - r_n.
\end{equation}
To proceed further, we equate the coefficients of $x^{-2}$ and $x^{-3}$ in $(qS_2^{\prime})$, which are long formulas. First equating the coefficients of $x^{-2}$ in $(qS_2^{\prime})$ produces
\begin{equation} \label{cxm2}
\begin{split}
 & \beta_n\left(q^{2n-1}-2\frac{q^{n-1}}{q^{\alpha}}(1-q)\left[\frac{t}{q}\right]\left(R_n+qR_{n-1} \right) \right)  = \\
& \qquad \left( \frac{1}{q^{2\alpha}q^{2n}}\frac{1-q^n}{1-q}\frac{1-q^{n+\alpha}}{1-q} + \left[ \frac{t}{q}\right]\frac{1}{q^{2\alpha}}\left( q^\alpha - \frac{2}{q^n}\right)(1-q^n)\right) + \frac{1}{q^{2\alpha}q^{n}}(2-q^n)(2-q^{n+\alpha})\left[ \frac{t}{q} \right]r_n\\
& \qquad \qquad \qquad + \left(\frac{1}{q^{2\alpha}q^n} - 2\frac{(1-q)^2}{q^{2\alpha}}\left[ \frac{t}{q} \right] \right)\left[ \frac{t}{q}\right]S_{n-1} + 2\frac{(1-q)^2}{q^{2\alpha}}\left[ \frac{t}{q} \right]^2 r_nS_{n-1}.
\end{split}
\end{equation}
Now equating the coefficients of $x^{-3}$ in $(qS_2^{\prime})$ gives
\begin{equation} \label{cxm3}
\begin{split}
& \beta_n\frac{1}{q^{\alpha}}q^{n-1}\left[ \frac{t}{q} \right]\left( R_n + qR_{n-1} \right)  = \\
& \qquad \qquad \frac{1}{q^{2\alpha}q^n}\left[ \frac{t}{q} \right]\left( \frac{1-q^n}{1-q}\right)  - \frac{1}{q^{2\alpha}}\frac{1}{q^n}\left(\frac{1-q^n+1-q^{\alpha+n}}{1-q} \right)\left[ \frac{t}{q}\right]r_n \\
&\qquad \qquad \qquad + \frac{1-q}{q^{2\alpha}}\left[ \frac{t}{q}\right]^2S_{n-1} -\frac{1-q}{q^{2\alpha}}\left[ \frac{t}{q}\right]^2r_nS_{n-1}.
\end{split}
\end{equation}
We now multiply (\ref{cxm3}) by $2(1-q)$ and add to (\ref{cxm2}) to obtain :
\begin{equation}
\beta_nq^{2n-1} = \frac{1}{q^{2\alpha}q^{2n}}\frac{1-q^n}{1-q}\frac{1-q^{n+\alpha}}{1-q} + \frac{1-q^n}{q^{\alpha}}\left[ \frac{t}{q }\right] + \frac{q^n}{q^{\alpha}}\left[ \frac{t}{q} \right]r_n + \frac{1}{q^{2\alpha}q^n}\left[ \frac{t}{q}\right]S_{n-1}.
\end{equation}
Meanwhile (\ref{cxm3}) can be written as
\begin{equation}
\begin{split}
\beta_nq^{n-1}\left( R_n + qR_{n-1} \right) &= \frac{1}{q^{\alpha}q^n}\left(\frac{1-q^n}{1-q} - \left( \frac{1-q^n}{1-q} + \frac{1-q^{n+\alpha}}{1-q}\right)r_n \right) \\
& \qquad \qquad \qquad \qquad \qquad + \frac{1-q}{q^{\alpha}}\left[ \frac{t}{q} \right] (1-r_n)S_{n-1}.
\end{split}
\end{equation}
In summary, so far we have the $5$ equations
\begin{equation} \label{1S}
q^{2n+\alpha} \alpha_n =\frac{(1-q^n)}{1-q}+\frac{1-q^{n+\alpha+1}}{q(1-q)}+ q^n\left[ \frac{t}{q} \right]\left(R_n+(1-q)S_{n-1}\right),
\end{equation}
\vskip 0.2cm
\noindent
\begin{equation} \label{2S}
r_{n+1} + r_n = -\alpha_nR_n +1,
\end{equation}
\vskip 0.6cm
\noindent
\begin{equation} \label{3S}
\beta_nq^{2n-1} = \frac{1}{q^{2\alpha}q^{2n}}\frac{1-q^n}{1-q}\frac{1-q^{n+\alpha}}{1-q}
+ \frac{1-q^n}{q^{\alpha}}\left[ \frac{t}{q }\right] + \frac{q^n}{q^{\alpha}}
\left[ \frac{t}{q} \right]r_n + \frac{1}{q^{2\alpha}q^n}\left[ \frac{t}{q}\right]S_{n-1},
\end{equation}
\vskip 0.6cm
\noindent
\begin{equation} \label{4S}
\begin{split}
\beta_nq^{n-1}\left( R_n + qR_{n-1} \right) &= \frac{1}{q^{\alpha}q^n}\left(\frac{1-q^n}{1-q} -
\left( \frac{1-q^n}{1-q} + \frac{1-q^{n+\alpha}}{1-q}\right)r_n \right) \\
& \qquad \qquad \qquad \qquad \qquad + \frac{1-q}{q^{\alpha}}\left[ \frac{t}{q} \right] (1-r_n)S_{n-1},
\end{split}
\end{equation}
\vskip 0.6cm
\noindent
\begin{equation} \label{5S}
\beta_nR_nR_{n-1} = r_n^2 - r_n.
\end{equation}

\begin{rem}
Substituting $t=0$ into (\ref{1S}) and (\ref{3S}) gives the recurrence coefficients for the $q-$Laguerre polynomials.
The recurrence coefficients for the Stieltjes-Wigert polynomials can be obtained from (\ref{1S}) and (\ref{3S}) by setting $\alpha=0$ and following the steps laid out in remark {\bf \ref{Askey-coef}}.
\end{rem}
\vskip 0.2cm
\noindent
Note that both equations (\ref{3S}) and (\ref{4S}) give an expression for $\beta_n$ in terms of the auxiliary quantities. Both of these equations are essential because they allow us to eliminate $\beta_n$ and obtain an expression for the sum $S_{n-1}$ in terms of $R_n$ and $r_n$ only.
The sum $S_{n-1}$ is given by the following lemma.
\begin{lemma}
If $S_{n} = \sum_{j=0}^n R_j$ then
\begin{equation} \label{sumR}
\begin{split}
& S_{n-1}\left[ \frac{t}{q}\right]\left( \frac{1}{q^{2\alpha}q^n}-\frac{q^n(1-q)(1-r_n)}{q^{\alpha}R_n} \right) = \\
& \qquad \qquad -\frac{1}{q^{2n+2\alpha}}\left( \frac{1-q^n}{1-q} \right)\left( \frac{1-q^{n+\alpha}}{1-q} \right) -\frac{q^n}{q^{\alpha}}\left[ \frac{t}{q}\right]r_n \\
& \qquad \qquad \qquad +\frac{1}{q^{\alpha}R_n}\left( \frac{1-q^n}{1-q}-\left(\frac{1-q^n}{1-q}+\frac{1-q^{n+\alpha}}{1-q} \right)r_n \right) \\
& \qquad \qquad \qquad \qquad -q^{2n}\frac{r_n^2-r_n}{R_n^2} - \frac{1-q^n}{q^{\alpha}}\left[ \frac{t}{q} \right].
\end{split}
\end{equation}
\end{lemma}
{\bf Proof}
\vskip .15cm
\noindent
First multiply (\ref{4S}) by $R_n$ and use (\ref{5S}) to eliminate $R_{n-1}$. Then substitute for $\beta_n$ from (\ref{4S}) into (\ref{3S}). This gives an expression for the sum $S_{n-1}$ in terms of $r_n$ and $R_n$ only.
\hfill $\square$
\vskip 0.2cm
\noindent
Note that this equation effectively eliminates the sum $S_{n-1}$ from equations (\ref{1S}) and (\ref{3S}) and consequently we see that $\alpha_n$ and $\beta_n$ are entirely determined by $r_{n}$ and $R_n$.

\section{Non-linear Difference equations}
\setcounter{equation}{0}
In this section we derive the coupled non linear difference equations given in theorem {\bf \ref{APl}}.
\vskip 0.2cm
\noindent
{\bf Proof of theorem \ref{APl}}
\vskip 0.2cm
\noindent
Eliminating $\alpha_n$ from (\ref{1S}) and (\ref{2S}) we obtain
\begin{equation} \label{diff1}
q^{2n+\alpha}(1-r_{n+1}-r_n) = \left(\frac{(1-q^n)}{1-q}+\frac{1-q^{n+\alpha+1}}{q(1-q)} + q^n\left[ \frac{t}{q} \right] \left(R_n+(1-q)S_{n-1}\right) \right)R_n.
\end{equation}
Eliminating $\beta_n$ from (\ref{4S}) and (\ref{5S}) we obtain
\begin{equation} \label{diff2}
\begin{split}
& q^{2n+\alpha-1} \left( r_n^2-r_n \right)\left( R_n + qR_{n-1} \right) = \\
& \qquad  \left(\frac{1-q^n}{1-q} - \left( \frac{1-q^n}{1-q} + \frac{1-q^{n+\alpha}}{1-q}\right)r_n + q^n(1-q)
\left[ \frac{t}{q} \right] (1-r_n)S_{n-1} \right)R_nR_{n-1}.
\end{split}
\end{equation}
We now replace $S_{n-1}$ in (\ref{diff1}) and (\ref{diff2}) by the expression given in (\ref{sumR}), and therefore obtain, respectively
\begin{equation} \label{rec1}
\begin{split}
& q^{n-\alpha}\left[ (1-r_{n+1})+\frac{1}{q}(1-r_n) \right]R_n - q^{3n}(1-q)(1-r_{n+1})(1-r_n) \\
& \qquad \qquad = 
\frac{t}{q^{2\alpha+1}}R_n^3 +\left(\frac{1}{q^{2\alpha+n+1}}\frac{1-q^{2n+\alpha+1}}{1-q}-(1-q)t\frac{q^n}{q^{\alpha+1}} \right)R_n^2 + q^{2n}R_n,
\end{split}
\end{equation}
and
\begin{equation} \label{rec2}
\begin{split}
& q^{2n-1}(1-q)r_n(1-r_n)^2 - \frac{R_n+qR_{n-1}}{q^{\alpha+1}}r_n(1-r_n) = \\
& \qquad \qquad R_nR_{n-1}\Bigg( \frac{1-q}{q^{\alpha+1}}t(1-r_n)(q^n(1-r_n)-1) + \frac{1-r_n}{q^{2n+2\alpha}}\left[ \frac{1-q^{2n+\alpha}}{1-q} \right] -\frac{1}{q^{2n+2\alpha}}\frac{1-q^{n+\alpha}}{1-q}  \Bigg).
\end{split}
\end{equation}
Note that (\ref{rec1}) is a first order difference equation in $r_n$ and a cubic in $R_n$, whilst (\ref{rec2}) is the other way around. \
These equations admit the respective factorizations,
\begin{equation} \label{factor1}
\begin{split}
& \left[R_n-q^{n+\alpha}(1-q)(q^{n+1}-q^{n+1}r_{n+1})\right]\left[R_n-q^{n+\alpha}(1-q)(q^{n}-q^nr_{n})\right] \\
&  \qquad \qquad \qquad  = -(1-q)q^{n}tR_n\left(R_n-q^{n+\alpha}(1-q)\right)\left(R_n-\frac{q^{n+\alpha+1}}{t(1-q)}\right)
\end{split}
\end{equation}
and
\begin{equation} \label{factor2}
\begin{split}
& q^nr_n\left[(q^n-q^nr_n)-\frac{R_n}{q^{n+\alpha}(1-q)}\right]\left[(q^n-q^nr_n)-\frac{R_{n-1}}{q^{n+\alpha-1}(1-q)}\right] \\
& \qquad \qquad = \frac{tR_nR_{n-1}}{q^{\alpha}}\left[(q^n-q^nr_n)-\frac{q}{t(1-q)^2}\right]\left[(q^n-q^nr_n)-1\right].
\end{split}
\end{equation}
Under the substitutions,
\begin{equation*}
x_n = \frac{q^{n+\alpha}(1-q)}{R_n}, \qquad \qquad y_n = q^n(1-r_n) \qquad \textrm{and} \qquad T = \frac{(1-q)^2}{q}t
\end{equation*}
(\ref{factor1}) and (\ref{factor2}) are the coupled equations in theorem {\bf \ref{APl}}.
\hfill $\square$

\section{Difference equation for $p(n)$ }
\setcounter{equation}{0}

We conclude this paper by stating without proof, a second order, non-linear difference equation, satisfied by
 the quantity $\varrho_n$, where
\begin{equation}
\varrho_n: = \frac{(1-q)^2}{q}p(n). 
\end{equation}
We are happy to provide the relevant PDF file upon request.
\vskip .2cm
\noindent
First let 
\begin{equation}
J_n: = \frac{1}{T}\left(1-q^{-2n-\alpha-1} + \frac{\varrho_n-q\varrho_{n+1}}{1-q} \right),
\end{equation}
\vskip 0.2cm
\noindent
and introduce the following variables, defined in terms of $J_n$ and $\varrho_n;$
\begin{equation}
F_n: = -q^{2n+\alpha}T\: J_n^2 + [q^{2n+\alpha}(\varrho_n+1)-1]J_n-q^n,
\end{equation}
\vskip 0.2cm
\begin{equation}
G_n: = q^{2n+\alpha}(T-\varrho_n)J_n^2+q^n(1-q^{n+\alpha})J_n,
\end{equation}
\vskip 0.2cm
\begin{equation}
H_n: = \frac{J_nJ_{n-1}q^n(1-q^{n+\alpha})}{J_n+J_{n-1}},
\end{equation}
\vskip 0.2cm
\begin{equation}
I_n: =q^n +  \frac{J_nJ_{n-1}(1-q^{2n+\alpha}(1+\varrho_n)) }{J_n+J_{n-1}}.
\end{equation}
After considerable computations, the second order, non-linear difference equation, satisfied by $\varrho_n$ reads,
\begin{equation}
\begin{split}
& (J_n+J_{n-1})^2G^2_n-(J_n+J_{n-1})\Big\{J_nJ_{n-1}q^n(1-q^{n+\alpha})(2G_n+F^2_n) \\
& \qquad +\left(q^n(J_n+J_{n-1})+J_nJ_{n-1}[1-q^{2n+\alpha}(1+\varrho_n)]\right)F_nG_n\Big\} \\
& \qquad \qquad +\Big\{q^n(J_n+J_{n-1})+J_nJ_{n-1}[1-q^{2n+\alpha}(1+\varrho_n)]\Big\}J_nJ_{n-1}q^n(1-q^{n+\alpha})F_n \\
& \qquad \qquad \qquad + \Big\{q^n(J_n+J_{n-1})+J_nJ_{n-1}[1-q^{2n+\alpha}(1+\varrho_n)] \Big\}^2G_n \\
& \qquad \qquad \qquad \qquad + (J_nJ_{n-1})^2q^{2n}(1-q^{n+\alpha})^2 = 0.
\end{split}
\end{equation}

\bibliographystyle{plain}

\end{document}